%%%% This is a standard AMS TeX source

\input amstex
\documentstyle{amsppt}

\NoBlackBoxes

\magnification=1200

\loadbold
\TagsOnRight
\mathsurround=1.2pt
\define\<{\snug}
\redefine\.{\,.}
\newdimen\theight
\define\lab#1{\vadjust{\setbox0=\hbox{\tenpoint%
            \quad${}^{\dsize\checkmark}{}^{\text{#1}}$}%
            \theight=\ht0
            \advance\theight by\dp0\advance\theight by \lineskip
            \kern-\theight\vbox to\theight{\rightline{\rlap{\box0}}%
            \vss}%
            }}%

\define\dddots{\mathinner{\mkern1mu\raise1pt\hbox{.}\mkern3mu
    \raise4pt\hbox{.}\mkern3mu\raise7pt\vbox{\kern7pt\hbox{.}}\mkern1mu}}

\redefine\>{\hskip1pt}
\redefine\le{\leqslant}
\redefine\ge{\geqslant}
\define\){)\hskip1pt}
\define\({\hskip1pt(}
\define\op{\operatorname}

\define\wt{\widetilde}
\define\g{\frak{g}}
\define\la{\lambda}
\define\e{\varepsilon}
\define\a{\alpha}
\define\be{\beta}
\define\lan{\langle}
\define\ran{\rangle}
\define\La{\Lambda}
\define\si{\sigma}
\define\th{\theta}
\define\dt{\delta}

\define\ph{\varphi}
\define\C{\Bbb C}
\define\mi{\>|\>}
\define\bcdot{\,\boldsymbol\cdot\,}

\topmatter
\title Shifted Schur functions II. \\
Binomial formula for characters of classical groups\\
and applications
\endtitle
\rightheadtext{Shifted Schur functions II}
\author Andrei~Okounkov and Grigori~Olshanski
\endauthor
\dedicatory To A.~A.~Kirillov on his 60\<th birthday
\enddedicatory
\address A. Okounkov: Dept. of Mathematics, University of
Chicago, 5734 South University Av., Chicago, IL 60637-1546.
E-mail address: okounkov\@math.uchicago.edu\newline
G. Olshanski: Institute for Problems of Information Transmission,
Bolshoy Karetny 19, 101447 Moscow GSP-4, Russia.
E-mail address: olsh\@ippi.ac.msk.su
\endaddress
\abstract Let $G$ be any of the complex classical groups
$GL(n)$, $SO(2n+1)$, $Sp(2n)$, $O(2n)$, let $\frak{g}$ denote the Lie
algebra of $G$, and let $Z(\frak{g})$ denote the subalgebra of
$G$-invariants in the universal enveloping algebra $U(\frak{g})$. We
derive a Taylor-type expansion for finite-dimensional characters of
$G$ (the binomial formula) and use it to specify a distinguished
linear basis in $Z(\frak{g})$. The eigenvalues of the basis elements
in highest weight $\frak{g}$-modules are certain shifted (or
factorial) analogs of Schur functions. We also study an associated
homogeneous basis in $I(\frak{g})$, the subalgebra of $G$-invariants
in the symmetric algebra $S(\frak{g})$. Finally, we show that the
both bases are related by a $G$-equivariant linear isomorphism
$\sigma\:I(\frak{g})\to Z(\frak{g})$, called the special
symmetrization.
\endabstract
\thanks The authors were supported by the Russian
Foundation for Basic Research (grant 95-01-00814).
The first author was supported by the NSF grant DMS 9304580.
\endthanks
\endtopmatter

\document

\head Introduction \endhead

The present paper is a continuation of our previous work \cite{OO1}
but can be read independently. Here we aim to transfer to the
orthogonal and symplectic groups a part of results of \cite{OO1}
connected with the general linear group.

Let $G$ be one of the groups $SO(2n+1,\Bbb C)$, $Sp(2n,\Bbb
C)$,\footnote{In the present Introduction, to simplify the
discussion, we exclude the even orthogonal group $SO(2n,\Bbb C)$.
However, after minor modifications, all our constructions hold for
this group as well.} let $\g$ be the Lie algebra of $G$, $U(\g)$ be
the universal enveloping algebra of $\g$, and $Z(\g)$ be the center
of $U(\g)$. The commutative algebra $Z(\g)$ is also known as the
algebra of biinvariant differential operators on the group $G$, or
the algebra of Laplace operators.

The irreducible finite-dimensional representations of the group $G$
are indexed by partitions $\la$ of length $\le n$; we denote them
by $V_\la$ and we view each $V_\la$ also as a $U(\g)$-module.

Our main object is a distinguished linear basis
$\{\Bbb{T}_\mu\}\subset Z(\g)$. Here the index $\mu$ ranges over
partitions of length $\le n$ and each basis element $\Bbb{T}_\mu$ can
be characterized, within a scalar multiple, as the unique element in
$Z(\g)$ of degree $2\>|\mu|$ satisfying the following {\it vanishing
condition\/}: the eigenvalue of $\Bbb{T}_\mu$ in $V_\la$ equals
$0$ if $|\la|\le|\mu|$ and $\la\ne\mu$ (moreover, it turns
out that the eigenvalue is $0$ if $\la_i<\mu_i$ for at least one
$i$).

For general $\la$, we denote the eigenvalue of $\Bbb{T}_\mu$ in
$V_\la$ by $t^*_\mu(\la)=t^*_\mu(\la_1,\dots,\la_n)$.
The functions $t^*_\mu$ are inhomogeneous polynomials with a certain
kind of symmetry; they can be identified with certain {\it
factorial\/} analogs of the Schur polynomials. Note that the top
homogeneous component $t_\mu$ of the inhomogeneous polynomial
$t^*_\mu$ coincides with a Schur polynomial in the squares of the
arguments:
$$
t_\mu(x_1,\dots,x_n)=s_\mu(x^2_1,\dots,x^2_n)\.
$$

The polynomials $t^*_\mu$ arise in the {\it binomial formula\/} for
the characters of $G$. Denote by
$\chi_\la=\chi_\la(z_1,\dots,z_n)$ the character of
$V_\la$; here $z_1,\dots,z_n$ are natural coordinates of a
maximal torus of $G$. Assume `additive' variables $x_1,\dots,x_n$ are
connected with the `multiplicative' variables $z_1,\dots,z_n$ by the
relation
$$
x_i=z^{1/2}-z^{-1/2},\quad\text{i.e.,}\quad x_i^2=z+z^{-1}-2\.
$$
The binomial formula is written as
$$
\frac{\chi_\la(z_1,\dots,z_n)} {\chi_\la(1,\dots,1)}=
\sum_\mu\frac{t^*_\mu(\la_1,\dots,\la_n\)t_\mu(x_1,\dots,x_n)}{c_\pm(n,\mu)},
$$
where $c_\pm(n,\mu)$ are simple normalization factors. This is a
Taylor-type expansion for the character near the unit element of the
group. \footnote{Note that the binomial formula for characters of
classical groups of type $B$, $C$, $D$ is a particular case of an
expansion of type $BC_n$ Jacobi polynomials, see Lassalle
\cite{Lass}.}

One more important property of the basis elements $\Bbb{T}_\mu$ is
existence of a relation between the bases corresponding to different
values of the parameter $n$. We call this the {\it coherence
property\/} of the canonical basis.

Then we study a counterpart of $\{\Bbb{T}_\mu\}$ for the algebra
$I(\g)$, the subalgebra of $G$-invariants in the symmetric algebra
$S(\g)$. Note that $I(\g)$ is a graded algebra which is canonically
isomorphic to $\operatorname{gr}Z(\g)$, the graded algebra
associated to the natural filtration in $Z(\g)$.

To the basis $\{\Bbb{T}_\mu\}$ corresponds a homogeneous basis
$\{T_\mu\}$ of $I(\g)$: each $T_\mu$ coincides with the leading term
of $\Bbb{T}_\mu$ (with respect to the identification
$I(\g)=\operatorname{gr}Z(\g)$). On the other hand, the basis
$\{T_\mu\}$ can be described without reference to $Z(\g)$, in
intrinsic terms of the $G$-module $S(\g)$. Note that the elements
$T_\mu$ admit an explicit expression as polynomials in natural
generators of the Lie algebra $\g$.

Finally, we show that there exists a $G$-equivariant linear
isomorphism
$$
\sigma\:S(\g)\to U(\g)
$$
which preserves leading terms and takes $T_\mu$ to $\Bbb{T}_\mu$ for
any $\mu$. We call $\sigma$ the {\it special symmetrization\/} for
the orthogonal/symplectic Lie algebras.

The special symmetrization $\sigma$ enters a family of `generalized
symmetrizations' studied in the paper \cite{O2}. The results of
\cite{O2} provide explicit combinatorial formulas both for $\sigma$
and its inverse $\sigma^{-1}$. Since we dispose with an explicit
expression for $T_\mu$, this yields a certain formula for the
elements $\Bbb{T}_\mu=\sigma(T_\mu)$.

We conclude this Introduction with a brief discussion of some related
works.

For the general linear group $GL(n,\Bbb C)$, the binomial theorem and
distinguished bases $\{\Bbb{S}_\mu\}\subset Z(\frak{gl}(n,\Bbb C))$
and $\{S_\mu\}\subset I(\frak{gl}(n,\Bbb C))$ with similar properties
were earlier considered in \cite{OO1}. The elements $\Bbb{S}_\mu$,
which are called {\it quantum immanants\/}, appear in a higher
version of the classical Capelli identity and admit a remarkable
explicit expression, see the papers \cite{Ok1, N,
Ok2}. The eigenvalues of the elements $\Bbb{S}_\mu$ in highest
weight modules are described by certain polynomials
$s^*_\mu(\la_1,\dots,\la_n)$, called in \cite{OO1} the {\it
shifted Schur polynomials\/}. The both families of polynomials,
$\{s^*_\mu\}$ and $\{t^*_\mu\}$, have similar properties.
\footnote{However, in contrast to the polynomials $s^*_\mu$, the
polynomials $t^*_\mu$ are not stable as $n\to\infty$.}

Let $\mu$ be of the form $(1^m)$ for the orthogonal group or of the
form $(m)$ for the symplectic group (\<$m=1,2,\dots$). Then the
corresponding elements $\Bbb{T}_\mu$ coincide with Laplace operators
considered by A.~Molev and M.~Nazarov in \cite{MN}. As shown in
\cite{MN}, these elements also occur in a Capelli-type identity and
admit an explicit expression in terms of the generators of the Lie
algebra. An open problem is to generalize the results of \cite{MN} to
arbitrary $\mu$. Note that the explicit formulas found in \cite{MN,
Ok1, N, Ok2} differ from the formulas obtained
via the special symmetrization.

Many of the results of \cite{OO1} and the present paper should have
counterparts for classical Lie superalgebras $\frak{gl}(p\mi q)$,
$\frak{q}(n)$, and $\frak{osp}(n\mi 2m)$. Some work in this direction
was made by A.~Borodin and N.~Rozhkovskaya \cite{BR}, A.~Molev
\cite{Mo}, V.~Ivanov and A.~Okounkov (see \cite{I}).

In the present work (as in \cite{OO1}), we used very elementary tools
such as explicit formulas for the (ordinary and factorial) Schur
polynomials or Weyl's character formula. We believe our approach can
be further developed by making use of a more fine technique. In
particular, there exist more involved versions of the binomial
formula, see \cite{OO2}, \cite{Ok3}, and \cite{Ok4}.

\head\S0.\enspace Main notation\endhead

Unless otherwise stated, the symbol $G(n)$ or $G$ will denote any of
the complex classical groups of rank $n$
$$
GL(n,\Bbb C),\quad Sp(2n,\Bbb C),\quad SO(2n+1,\Bbb C),\quad SO(2n,\Bbb C),
$$
which constitute the series A, C, B, D, respectively. By
$\g(n)$ or $\g$ we denote the corresponding complex Lie
algebras. In the case of the series D we also consider the
nonconnected groups $G'=G'(n)=O(2n,\Bbb C)$.

By $U(\g)$ we denote the universal enveloping algebra of $\g$.
For the series A, C, B, we denote by $Z(\g)$ the
center of $U(\g)$ (which coincides with the subalgebra of $G$-invariants
in $U(\g)$); for the series D the same symbol
will denote the subalgebra of $G'$-invariants in $U(\g)$
(which is a proper subalgebra of the center).

Similarly, we define $I(\g)\subset S(\g)$ as the subalgebra
of $G$-invariants (or $G'$-invariants, for the series D) in
the symmetric algebra of $\g$.

$V_\la$ or $V_{\la\mi n}$ is the irreducible finite-dimensional
complex-analytic representation of $G$ with highest weight $\la$; it
is also viewed as a $U(\g)$-module.

$\C\>[G]$ is the space of regular functions on $G$, where $G$ is
viewed as an algebraic group over $\C$; in other words, $\C\>[G]$
is the linear span of matrix elements of all representations
$V_{\la}$.

By $\mu$ we always denote a partition of length $l(\mu)\le n$; it is
also viewed as a Young diagram; $|\mu|$ denotes the number of boxes of
the diagram $\mu$.

\head\S1.\enspace Binomial formula\endhead

Here we derive a multidimensional analog of the binomial formula
$$
(1+x)^k=\sum_{m=0}^k\frac1{m!}\,k(k-1)\cdots(k-m+1\)x^m,
$$
where the powers of a variable will be replaced by characters of $G$.

Let us agree about the choice of a Borel subgroup $B\subset G$ and a
maximal torus $H\subset B$ (thus positive roots and dominant weights will
be specified):

$\bullet$ For the series A, we take as $B$ the subgroup of
upper triangular matrices and as $H$ the subgroup of diagonal
matrices,
$$
H=\{\op{diag}\(z_1,\dots,z_n)\},\qquad z_1,\dots,z_n\in\C^*.
$$

$\bullet$ For the series C--B--D we identity $G$ with a subgroup in
$GL(N,\C)$, where $N=2n$ or $N=2n+1$,
$$
G=\{g\in GL(N,\C)\mid g'Mg=M\},\tag1.1
$$
where $M$ stands for the following symmetric (case B--D) or
antisymmetric (case C) matrix of order $N$:
$$
M=\left[\matrix 0&\matrix&&1\\
&\dddots&\\
1&&\endmatrix\\
\matrix&&\pm1\\
&\dddots&\\
\pm1&&\endmatrix&0\endmatrix\right].\tag1.2
$$

Then as $B$ and $H$ we take the subgroups of upper triangular or
diagonal matrices in $G$. Thus
$$
\alignat2
&H=\{\op{diag}\(z_1,\dots,z_n,z_n^{-1},\dots,z_1^{-1})\},&\qquad
&\text{case C--D}\\
&H=\{\op{diag}\(z_1,\dots,z_n,1,z_n^{-1},\dots,z_1^{-1})\},&\qquad
&\text{case B}\.
\endalignat
$$

For all the series, the weights of $H$ are of the form
$$
z^{\la}=z_1^{\la_1}\dots z_n^{\la_n},\quad\text{where }\,
\la=(\la_1,\dots,\la_n)\in\Bbb Z^n,
$$
while the dominant weights are
distinguished by the supplementary conditions
$$
\alignat2
&\la_1\ge\cdots\ge\la_n,&\qquad&\text{case A},\\
&\la_1\ge\cdots\ge\la_n\ge0,&\qquad&\text{case C--B},\\
&\la_1\ge\cdots\ge\la_{n-1}\ge|\la_n|,&\qquad&\text{case D}\.
\endalignat
$$

The dominant weights for the series A are called the {\it
signatures\/} and those for the series C--B are called the {\it
positive signatures\/} (\<$=$partitions of length $\le n$).

The irreducible character of $G$, indexed by a dominant weight $\la$
(i.e., the character of $V_{\la}=V_{\la|n}$) will be denoted by the
symbol
$$
\chi_{\la}=\chi_{\la}^{gl(n)},\quad\chi_{\la}^{sp(2n)},\quad
\chi_{\la}^{so(2n+1)}\quad\text{or }\,\chi_{\la}^{so(2n)}.
$$
For the series D, it is more convenient to deal with the {\it
reducible\/} character
$$
\chi_{\la}^{o(2n)}=\chi_{\la_1,\dots,\la_n}^{so(2n)}+\chi_{\la_1,\dots,
\la_{n-1},-\la_n}^{so(2n)},\quad\text{where }\,\la_1\ge\cdots\ge\la_n\ge0\.
$$
When $\la_n=0$, this is simply $2\chi_{\la}^{so(2n)}$, and when
$\la_n>0$, this coincides with the restriction to the group
$G=SO(2n,\C)$ of an irreducible character of the group
$G'=O(2n,\C)$.

Thus, we are working with the characters
$$
\chi_{\la}^{gl(n)},\quad\chi_{\la}^{sp(2n)},\quad\chi_{\la}^{so(2n+1)},\quad
\chi_{\la}^{o(2n)},
$$
where $\la$ is either a signature (case A) or a positive signature
(case C--B--D).

To each $\la$ we assign another weight $l=(l_1,\dots,l_n)$ defined as
$$
\alignat2
&l=(\la_1+n-1,\,\la_2+n-2,\dots,\la_n),&\qquad&\text{case A},\tag1.3\\
&l=(\la_1+n-1+\e,\,\la_2+n-2+\e,\dots,\la_n+\e),&\qquad
&\text{case C--B--D},\tag1.4
\endalignat
$$
where
$$
\e=\cases1,&\text{case C},\\
1/2,&\text{case B},\\
0,&\text{case D}\.\endcases\tag1.5
$$
For the series C--B--D, $l$ just equals $\la$ plus the half-sum of
positive roots, and for the series A, the same holds after
restriction to $SL(n,\C)\subset GL(n,\C)$.

Unless otherwise stated, we shall consider characters of $G$ as
functions on the torus $H$, i.e., as functions of the coordinates
$z_1,\dots,z_n$. The following formulas follow from general Weyl's
character formula (the determinants below are of order~$n$):
$$
\align
\chi_{\la}^{gl(n)}(z_1,\dots,z_n)&=\frac{\det\>[z_j^{l_i}]}{\Pi_{i<j}
(z_i-z_j)},\\
\chi_{\la}^{sp(2n)}(z_1,\dots,z_n)&=\frac{\det\>[(z_j^{l_i}
-z_j^{-l_i})/(z_j-z_j^{-1})]}{\Pi_{i<j}(z_i+z_i^{-1}-z_j-z_j^{-1})},\\
\chi_{\la}^{so(2n+1)}(z_1,\dots,z_n)&=\frac{\det\>[(z_j^{l_i}
-z_j^{-l_i})/(z_j^{1/2}-z_j^{-1/2})]}{\Pi_{i<j}(z_i+z_i^{-1}-z_j-z_j^{-1})},\\
\chi_{\la}^{o(2n)}(z_1,\dots,z_n)&=\frac{\det\>[z_j^{l_i}+z_j^{-l_i}]}
{\Pi_{i<j}(z_i+z_i^{-1}-z_j-z_j^{-1})}\.
\endalign
$$

Note that $\chi_{\la}^{gl(n)}(z_p,\dots,z_n)=s_{\la}(z_1,\dots,z_n)$,
the Schur polynomial in $n$ variable, if $\la$ is a partition. The
characters $\chi_{\la}^{gl(n)}$ are symmetric Laurent polynomials in
$z_1,\dots,z_n$, while the other characters are symmetric polynomial
functions of the variables $z_i+z_i^{-1}$, $i=1,\dots,n$. Note that
the last claim fails when $\chi_{\la}^{o(2n)}$ is replaced by
$\chi_{\la}^{so(2n)}$; this is the reason why we deal with the
characters $\chi_{\la}^{o(2n)}$.

Let $a=(a_1,a_2,\dots)$ be an arbitrary number sequence.
The {\it generalized factorial powers\/} of a variable $x$ are defined as
$$
(x\mi a)^k=\cases (x-a_1)\cdots(x-a_k),&k\ge1,\\
1,&k=0\.\endcases
$$
When $a\equiv0$ these are the ordinary powers, when $a=(0,1,2,\dots)$
these are the {\it falling factorial powers}, and when $a=(0,-1,-2,\dots)$
these are the {\it raising factorial powers\/} (or the Pohgammer symbol).

Further, the {\it generalized factorial Schur polynomial\/} in $n$
variables, indexed by a partition $\mu$ with $l(\mu)\le n$, is defined as
$$
s_\mu(z_1,\dots,x_n\mi a)
=\frac{\det\>[(x_j\mi a)^{\mu_i+n-i}]}{\det\>[(x_j\mi a)^{n-i}]}
=\frac{\det\>[(x_j\mi a)^{\mu_i+n-i}]}{\prod_{i<j}(x_i-x_j)}\.
$$

When $a\equiv0$ this turns into the standard formula for the
ordinary Schur polynomial $s_\mu(x_1,\dots,x_n)$, and in the general case
$s_\mu(x_1,\dots,x_n\mi a)$ is an {\it inhomogeneous\/} symmetric
polynomial of degree $|\mu|=\mu_1+\cdots+\mu_n$ and with the top
homogeneous component equal to $s_\mu(x_1,\dots,x_n)$:
$$
s_\mu(x_1,\dots,x_n\mi a)=s_\mu(x_1,\dots,x_n)+\text{lower terms}\.
$$
It follows that the polynomials $s_\mu(x_1,\dots,x_n\mi a)$ form a
basis in the algebra of symmetric polynomials in $n$ variables

The polynomials $s_\mu(x_1,\dots,x_n\mi 0,1,\dots)$ were introduced
by Biedenharn and Louck and called the {\it factorial Schur
polynomials}, see \cite{BL1, BL2}. The general definition is due
to Macdonald, see \cite{M1} and \cite{M2, I, \S3, ex.~20}
(note that Macdonald
uses the notation $(x\mi a)^k=(x+a_1)\cdots(x+a_k)$).

\proclaim{Theorem 1.1 \rm(Binomial formula for $GL(n)$)} Assume
$$
z_1=1+x_1,\quad\dots,\quad z_n=1+x_n\.
$$
Then
$$
\frac{\chi_\la^{gl(n)}(z_1,\dots,z_n)}{\chi_{\la}^{gl(n)}(1,\dots,1)}=
\sum_{\mu}\frac{s_{\mu}(l_1,\dots,l_n\mi0,1,2,\dots\)
s_{\mu}(x_1,\dots,x_n)}{c(n,\mu)}\.\tag1.6
$$
Here $\la=(\la_1,\dots,\la_n)$ is an arbitrary signature,
$l_i=\la_i+n-i$ \rom(\<$i=1,\dots,n$\<\rom),
$\mu$ ranges over partitions of
length $\le n$, and
$$
c(n,\mu)=\prod_{(i,j)\in\mu}(n+j-i)=\prod_{i=1}^n\frac{(\mu_i+n-i)!}{(n-i)!},
\tag1.7
$$
where the first product is taken over all boxes $(i,j)$ of the Young
diagram representing the partition $\mu$.
\endproclaim

\remark{Comment} This version of the binomial formula was proposed by the
authors, see \cite{OO1, Theorem 5.1}. It is equivalent to the
expansion
$$
s_{\la}(1+x_1,\dots,1+x_n)=\sum_\mu d_{\la\mu}s_\mu(x_1,\dots,x_n),
$$
where the coefficients $d_{\la\mu}$ are given by
$$
d_{\la\mu}=\det\left[\pmatrix\la_i+n-i\\
\mu_j+n-j\endpmatrix\right]_{1\le i,j\le n},
$$
see Lascoux \cite{Lasc}, Macdonald \cite{M2, I, \S3, ex.~10}.
\endremark

Although the latter formula looks simpler than (1.6), formula (1.6)
has a number of advantages, as it is explained in our paper
\cite{OO1}. E.g., an important idea is to consider
$l_1,\dots,l_n$ as variables rather than parameters, and then formula
(1.6) makes in evidence the fact that the roles of $l_1,\dots,l_n$
and of $x_1,\dots,x_n$ are almost symmetric.

\proclaim{Theorem 1.2 \rm(Binomial formula for the series C--B--D)} Assume
the variables $z_i$ and $x_i$ are subject to the relations
$$
x_i^2=z_i+z_i^{-1}-2,\quad\text{i.e.},\quad
x_i=\pm(z_i^{1/2}-z_i^{-1/2}),\qquad i=1,\dots,n\.
$$
\endproclaim

Let $G$ be any of the groups of the series C, B, D and let us
abbreviate $\chi_{\la}$ for $\chi_{\la}^{sp(2n)}$,
$\chi_{\la}^{so(2n+1)}$ or $\chi_{\la}^{o(2n)}$, where $\la$ is a
partition of length $\le n$. Then we have
$$
\frac{\chi_{\la}(z_1,\dots,z_n)}{\chi_{\la}(1,\dots,1)}
=\sum_{\mu}\frac{s_{\mu}(l_1^2,\dots,l_n^2\mi\e^2,(\e+1)^2,\dots\)
s_{\mu}(x_1^2,\dots,x_n^2)}{c_{\pm}(n,\mu)}\.\tag1.8
$$
Here $l_i=\la_i+n-i+\e$, $\e$ is defined by (1.5), $\mu$ ranges over
partitions of length $\le n$, and
$$
c_{\pm}(n,\mu)=\prod_{(i,j)\in\mu}4\(n+j-i)(n\pm1/2+j-i),\tag1.9
$$
where the sign ``$+$'' is taken for the series C, B, while the
sign ``$-$'' is taken for the series D.

Formula (1.8) is a particular case of the expansion of a type $BC_n$
Jacobi polynomial with ``Jack parameter'' $\a=1$ into a series of
Schur polynomials. This expansion probably is well-known to experts
(see, e.g., Lassalle \cite{Lass, Th\'eor\`eme~9}). For sake of
completeness we present a detailed proof.

\demo{Proof} Consider the classical Jacobi polynomials $P_k^{(\a,\be)}$
and put
$$
\wt P_k^{(\a,\be)}=P_k^{(\a,\be)}/P_k^{(\a,\be)}(1)\.
$$
Let $T_k$ and $U_k$ be the Chebyshev polynomials of the first and
second kind, respectively. Then for any $k=0,1,\dots$
$$
\alignat2
\frac{z^{k+1}-z^{-k-1}}{z-z^{-1}}
&=U_k\bigg(\frac{z+z^{-1}}2\bigg)&
&\sim\wt{P}^{(1/2,1/2)}_k\bigg(\frac{z+z^{-1}}2\bigg)\\
\frac{z^{k+1/2}-z^{-k-1/2}}{z^{1/2}-z^{-1/2}}&&
&\sim\wt{P}^{(1/2,-1/2)}_k\bigg(\frac{z+z^{-1}}2\bigg)\\
z^k+z^{-k}&=2T_k\bigg(\frac{z+z^{-1}}2\bigg)&
&\sim\wt{P}^{(-1/2,-1/2)}_k\bigg(\frac{z+z^{-1}}2\bigg),
\endalignat
$$
where the symbol ``$\sim$'' means equality within a number factor
depending on $k$ but not depending on $z$ (see Szeg\"o \cite{S,
(4.1.7) and (4.1.8)}).

Substituting these expressions into the formulas for the characters
$\chi_\la$ given above we obtain
$$
\chi_{\la}(z_1,\dots,z_n)\sim\frac{\det\>
[\wt{P}_{\la_i+n-i}^{(\a,\be)}(1+t_j)]_{1\le i,j\le
n}}{\prod_{i<j}(t_i-t_j)},\tag1.10
$$
where
$$
t_i=\frac{z_i+z_i^{-1}}2-1,\quad i=1,\dots,n,\qquad
(\a,\be)=\cases(1/2,1/2),&\text{case C},\\
(1/2,-1/2),&\text{case B},\\
(-1/2,-1/2),&\text{case D},\endcases
$$
and the symbol ``$\sim$'' means equality within a factor not
depending on $t_1,\dots,t_n$.
\enddemo

Recall a well-known identity: for any power series
$$
f_i(t)=\sum_{m=0}^{\infty}a_m^{(i)}t^m,\qquad i=1,\dots,n,
$$
we have
$$
\frac{\det\>[f_i(t_j)]_{1\le i,j\le n}}
{\prod_{i<j}(t_i-t_j)}=\sum_{\mu}\det\>[a_{\mu_j+n-j}^{(i)}]\>s_{\mu}
(t_1,\dots,t_n),\tag1.11
$$
where $\mu$ ranges over partitions of length $\le n$ (see, e.g., Hua
\cite{Hua, Theorem 1.2.1}).

We shall apply this identity to
$$
f_i(t)=\wt{P}_{\la_i+n-i}^{(\a,\be)}(1+t),\qquad i=1,\dots,n,
$$
and we shall use the well-known expansion
$$
\wt P_k^{(\a,\be)}(1+t)=\sum_{m\ge
0}\frac{k(k-1)\cdots(k-m+1)(k+\a+\be+1)\cdots(k+\a+\be+m)}
{2^mm!\,(\a+1)\cdots(\a+m)}\,t^m
$$
(see, e.g., Szeg\"o \cite{S, (4.21.2)}).

Remark that in our situation $\a+\be+1=2\e$, so that
$$
(k-i)(k+\a+\be+1+i)=(k+\e)^2-(\e+i)^2,\qquad i=0,1,\dots\.
$$
Further, remark that $\a=\pm1/2$ according to the assumption on the
sign made in the statement of the theorem. It follows that the above
expansion can be rewritten as
$$
\wt P_k^{(\a,\be)}(1+t)=\sum_{m\ge0}
\frac{((k+\e)^2\mi\e^2,(\e+1)^2,\dots)^mt^m}{2^mm!\,(\pm1/2+1)
\cdots(\pm1/2+m)}\.\tag1.12
$$

Formulas (1.10), (1.11), and (1.12) imply that
$$
\split
&\chi_{\la}(z_1,\dots,z_n)\\
&\qquad\sim\sum_{\mu}\frac{\det\>[(l_i^2\mi
\e^2,(\e+1)^2,\dots)^{\mu_j+n-j}]\>s_{\mu}(t_1,\dots,t_n)}
{\prod_{i=1}^n(2^{\mu_i+n-i}(\mu_i+n-i)!\,
(\pm1/2+1)\cdots(\pm1/2+\mu_i+n-i))}\.
\endsplit
$$
The normalizing factor in this formula is determined by the condition
that the constant term (i.e., the coefficient of
$s_\varnothing(t_1,\dots,t_n)\equiv1$) should be equal to 1. After simple
transformations, using the relations
$$
\gather
s_{\mu}(t_1,\dots,t_n)=2^{-|\mu|}s_{\mu}(x_1^2,\dots,x_n^2),\\
\frac{\det\>[(l_i^2\mi\e^2,(\e+1)^2,\dots)^{\mu_j+n-j}]}{\det\>[(l_i^2\mi
\e^2,(\e+1)^2,\dots)^{n-j}]}
=s_{\mu}(l_1^2,\dots,l_n^2\mi\e^2,(\e+1)^2,\dots),
\endgather
$$
we obtain the desired formula.\qed

\head\S2.\enspace A distinguished basis in $Z(\g)$\endhead

The aim of this section is to construct and characterize a linear
basis in $Z(\g)$. Our construction is suggested by the binomial
formula of \S1.

Let $\Cal O_e(G)$ be the algebra of germs of holomorphic functions at the
unity $e$ of the group $G$ and let $\Cal{M}_e(G)$ be the maximal ideal of
$\Cal O_e(G)$ formed by germs vanishing at $e$. We may identity $U(\g)$
with the space of those linear functionals on $\Cal O_e(G)$ that
vanish on $\Cal{M}_e(G)^m$, where $m$ is large enough. This defines a
nondegenerate pairing $\lan\bcdot,\bcdot\ran$ between $U(\g)$ and
$\C\>[G]$ (see \S0 for the definition of $\C\>[G]$). If $V$ is a
finite-dimensional $G$-module (also viewed as a $U(\g)$-module),
$V^*$ is the dual module, $\xi\in V$ and $\eta\in V^*$
are arbitrary vectors, and $f_{\xi\eta}(g)=\eta(g\xi)$ is the
corresponding matrix coefficient, then we have
$$
\lan X,f_{\xi\eta}\ran=\eta(X\xi)\quad\text{for any }\,X\in U(\g)\.
$$

Let $I(G)\subset\C\>[G]$ be the subalgebra of invariants of the
group $G$ (for series A, C, B) or the group $G'$ (for series
D) with respect to its action on $G$ by conjugations; note that
this definition is parallel to that of $Z(\g)\subset U(\g)$, see \S0.

Remark that both $U(\g)$ and $\C\>[G]$ are semisimple modules
over $G$ or $G'$. It follows that there are canonical projections
$$
\align
&\#\:U(\g)\to Z(\g),\tag2.1\\
&\#\:\C\>[G]\to I(G);\tag2.2\endalign
$$
these projections also can be defined as averaging over a compact
form of the group $G$ or $G'$. The pairing $\lan\bcdot,\bcdot\ran$ over
$U(\g)$ and $\C\>[G]$ is invariant over $G$ or $G'$ and so
defines a nondegenerate pairing between $Z(\g)$ and $I(G)$.

\proclaim{Lemma-Definition 2.1} Let $\mu$ range over the set of partitions
of length $\le n$. There exist central elements
$$
\alignat3
\Bbb S_{\mu}&\in Z(\g),&\quad\op{deg}\Bbb S_{\mu}&\le|\mu|&\qquad&
(\text{series A}),\\
\Bbb T_{\mu}&\in Z(\g),&\quad\op{deg}\Bbb T_{\mu}&\le2\>|\mu|&\qquad&
(\text{series C, B, D}),
\endalignat
$$
uniquely specified by the condition
$$
\align
&(f^\#\bigm|_H)(z_1,\dots,z_n)=\sum_{\mu}\frac{\lan\Bbb S_{\mu},f\ran\>
s_{\mu}(x_1,\dots,x_n)}{c(n,\mu)},\\
&(f^\#\bigm|_H)(z_1,\dots,z_n)=\sum_{\mu}\frac{\lan\Bbb T_{\mu},f\ran\>
s_{\mu}(x_1^2,\dots,x_n^2)}{c_{\pm}(n,\mu)},\endalign
$$
where $H\subset G$ is the diagonal torus defined in\/ \rom{\S1},
$f\in\C\>[G]$ stands for an arbitrary test function\/\rom;
$f^\#\in I(G)$ is the
image of $f$ under the projection\/ \rom{(2.2)}, and the variables
$x_1,\dots,x_n$ and the normalizing factors $c(n,\mu)$,
$c_{\pm}(n,\mu)$ are given by\/\rom{(1.6)}, \rom{(1.9)}.
\endproclaim

\demo{Proof} The function $f^\#\bigm|_H$ is a symmetric Laurent polynomial in
$z_1,\dots,z_n$; in the case of the series C, B, D it is also
invariant under transformations $z_i\mapsto z_i^{-1}$. It follows
that $f^\#\bigm|_H$ can be expanded into a series of polynomials
$s_{\mu}(x_1,\dots,x_n)$ or $s_{\mu}(x_1^2,\dots,x_n^2)$, where
$x_i=z_i-1$ (for series A) or $x_i^2=z_i+z_i^{-1}-2$ (for series
C, B, D). Given $\mu$, the map assigning to $f$ the $\mu$\<th
coefficient in that expansion is a linear functional on $\C\>[G]$
which only depends on a finite jet of $f$ at the point $e\in G$ (the
order of the jet is equal to $|\mu|$ or $2\>|\mu|$). Hence this
functional is an element of the algebra $U(\g)$ of degree
$\le|\mu|$ or $\le2\>|\mu|$. This proves existence of the elements
$\Bbb S_{\mu}$ or $\Bbb T_{\mu}$. As will be clear in what follows,
we have exact equalities
$$
\op{deg}\Bbb S_{\mu}=|\mu|,\qquad\op{deg}\Bbb T_{\mu}=2\>|\mu|\.
$$
Finally, the elements $\Bbb S_{\mu}$, $\Bbb T_{\mu}$ belong to
$Z(\g)$, because their values at any test function $f\in\C\>[G]$
depend on $f^\#$ only.\qed
\enddemo

Recall that $V_{\la}$ denotes the irreducible $G$-module with highest
weight $\la=(\la_1,\dots,\la_n)$, which is equally viewed as a
$U(\g)$-module. Any central element $X\in U(\g)$ acts in
$V_\la$ as a scalar operator $\op{const}\cdot1$;
the corresponding constant is
called the {\it eigenvalue of the central element $X$ in $V_{\la}$}.

\proclaim{Theorem 2.2} Let $G=GL(n,\C)$. The elements $\Bbb S_{\mu}$
defined in Lemma\/ \rom{2.1} form a basis in $Z(\g)$ and are
characterized\/ \rom(within a scalar multiple\/\rom)
by the following properties\/\rom:
\roster
\item"(\rom{i})" $\op{deg}\Bbb S_{\mu}\le|\mu|$\rom;
\item"(\rom{ii})" the eigenvalue of
$\Bbb S_{\mu}$ in a module $V_{\la}$ equals $0$ for all partitions
$\la$ such that $l(\la)\le n$, $|\la|\le|\mu|$, $\la\ne\mu$\rom;
\item"(\rom{iii})" the eigenvalue of $\Bbb S_{\mu}$ in $V_{\mu}$ is nonzero.
\endroster
\endproclaim

\proclaim{Theorem 2.3} Let $G$ be a classical group of type $C,B,D$. The
elements $\Bbb T_{\mu}$ defined in Theorem\/ \rom{2.1} form a basis in
$Z(\g)$ and are characterized\/ \rom(within a scalar multiple\/\rom) by the
following properties\/\rom:
\roster
\item"(\rom{i})" $\op{deg}\Bbb T_{\mu}\le2\>|\mu|$\rom;
\item"(\rom{ii})" the eigenvalue of $\Bbb T_{\mu}$ in a module $V_{\la}$
equals $0$ for all partitions $\la$ such that $l(\la)\le n$,
$|\la|\le|\mu|$, $\la\ne\mu$\rom;
\item"(\rom{iii})" the eigenvalue of $\Bbb T_{\mu}$ in $V_{\mu}$ is nonzero.
\endroster
\endproclaim

\remark{Comment} Note a difference between these two claims: in Theorem
2.3 we used all irreducible $G$-modules\footnote{The fact that for the
series D there exist dominant highest weights
$\la=(\la_1,\dots,\la_n)$ with $\la_n<0$ is not relevant here,
because for any element of $Z(o(2n))$, its eigenvalue does not change
under the transformation $\la_n\mapsto-\la_n$.} to characterize our
basis, while in Theorem 2.2 we used only a part of the modules
(namely, the polynomial ones). This implies that for the series A,
our basis $\{\Bbb S_{\mu}\}$ is not invariant under the outer
automorphism of $\frak{gl}(n,\C)$.
\endremark

To prove Theorems 2.2--2.3 we shall restate them as purely
combinatorial claims. We shall need the Harish--Chandra homomorphism.
Recall (see, e.g., Dixmier \cite{D, 7.4}) that it is defined for any
reductive Lie algebra $\g$ over $\C$. Fix a triangular
decomposition $\g=\goth n_-\oplus\goth{h}\oplus\goth n_+$ where
$\goth{h}$ is a Cartan subalgebra and $\goth n_+$ and $\goth n_-$ are
spanned by positive and negative root vectors. The Harish--Chandra
homomorphism maps $U(\g)^{\goth{h}}$, the centralizer of $\goth{h}$
in $U(\g)$, onto $U(\goth{h})=S(\goth{h})$. Let us identity
$S(\goth{h})$ with $\C\>[\goth{h}^*]$, the algebra of polynomial
functions on the dual space $\goth{h}^*$. It is known that the restriction
of the Harish--Chandra homomorphism to the center of
$U(\g)$ is injective and its image in $\C\>[\goth{h}^*]$
consists of those polynomials $f(\la)$ on $\goth{h}^*$ that are
invariant under transformations
$$
\la\mapsto w(\la+\rho)-\rho,\qquad\la\in\goth{h}^*,\tag2.1
$$
where $w$ ranges over the Weyl group $W=W(\g,\goth{h})$ and
$\rho$ denotes the half-sum of positive roots. Moreover, if $X$ is a
central element of $U(\g)$ and $f_X$ is the corresponding
invariant polynomial then the value of $f_X$ at a point $\la\in\goth{h}^*$
coincides with the eigenvalue of $X$ in the irreducible highest
weight module over $\g$ with highest weight $\la$. Although this
is true for all $\la\in\goth{h}^*$ we shall deal with dominant weights
$\la$ only; note that any polynomial $f(\la)$ is uniquely determined
by its values on the set of dominant weights. Finally, note that
$\op{deg}X=\op{deg}f_X$.

After this digression let us return to the classical Lie algebras.
Examine first the series A. Then the group $W$ coincides with the
symmetric group $S(n)$ and the weight $\rho$ satisfies the property
$\rho_i-\rho_{i+1}=1$, $1\le i\le n-1$. It follows that a polynomial
$f(\la_1,\dots,\la_n)$ is invariant under the action (2.1) if and only
if it satisfies the following symmetry property:
$$
f(\la_1,\dots,\la_n)=f(\la_1,\dots,\la_{i-1},\la_{i+1}-1,
\la_i+1,\la_{i+2},\dots,\la_n)\tag2.2
$$
for $i=1,\dots,n-1$. In our paper \cite{OO1} such polynomials were called
{\it shifted symmetric\/} ones and the algebra of shifted symmetric
polynomials in $n$ variables was denoted by $\La^*(n)$. Clearly,
$f(\la_1,\dots,\la_n)$ is shifted symmetric if and only if it is
symmetric in the variables $l_1,\dots,l_n$ that are related to
$\la_1,\dots,\la_n$ by (1.3). Note that the top homogeneous component
of a shifted symmetric polynomial is a symmetric polynomial. The
Harish--Chandra homomorphism induces an algebra isomorphism
$Z(\frak{gl}(n))\to\La^*(n)$.

Next assume $G$ is one of the classical groups of type C, B, D. Then
the weight $\rho$ has the form
$$
\rho=\rho_{\e}=(n-1+\e,\,n-2+\e,\dots,\e),\quad\text{where }\,\e=1,1/2,0,
$$
and the group $W$ is either the hyperoctahedral group
$H(n)=S(n)\ltimes\{\pm1\}^n$ (C--B case) or its subgroup
$H'(n)\subset H(n)$ of index $2$ (D case); recall that
the action of the generators of the subgroup $\{\pm1\}^n$ is of the
form $\la_i\mapsto-\la_i$. However, for the series D, we have
replaced the center of $U(o(2n,\C))$ by its proper subalgebra
$Z(o(2n,\C))$, which is just equivalent to replacement of $H'(n)$
by the whole group $H(n)$. This allows us to describe the algebra
$Z(\g)$ in a uniform way for the three series C, B, D.

Denote by $\roman{M}^*(n)=\roman{M}^*_{\e}(n)$ the algebra of polynomials in
$\la_1,\dots,\la_n$ that can be written in the form
$$
f(\la_1,\dots,\la_n)=g(l_1^2,\dots,l_n^2),\qquad l=\la+\rho_{\e},
$$
where $g$ is an arbitrary symmetric polynomial. Then we obtain an
algebra isomorphism $Z(\g)\to \roman{M}^*(n)$ induced by the
Harish--Chandra homomorphism.

\proclaim{Lemma-Definition 2.4} For the series A, put
$$
s^*_{\mu}(\la_1,\dots,\la_n)=s_{\mu}(l_1,\dots,l_n\mi0,1,2,\dots)
$$
and for the series C, B, D, put
$$
t^*_{\mu}(\la_1,\dots,\la_n)=s_{\mu}(l_1^2,\dots,l_n^2\mi\e^2,(\e+1)^2,
\dots),
$$
where $\e=1,1/2,0$, respectively. The Harish--Chandra homomorphism
sends $\Bbb S_{\mu}$ to $s^*_{\mu}$ and $\Bbb T_{\mu}$ to $t^*_{\mu}$.
\endproclaim

\demo{Proof} This follows at once from the above discussion if one
compares the definition of the elements $\Bbb S_{\mu}$, $\Bbb
T_{\mu}$ (Lemma-Definition 2.1), the claim of Theorems 1.1--1.2, and
the above definition of $s^*_{\mu}$ and $t^*_{\mu}$.\qed
\enddemo

Since the map $Z(\frak{gl}(n))\to\La^*(n)$ preserves degree and since
$\op{deg}s^*_{\mu}=|\mu|$, we have $\op{deg}\Bbb S_{\mu}=|\mu|$.
Similarly, since the top homogeneous component of
$t^*_{\mu}(\la_1,\dots,\la_n)$ coincides with
$s_{\mu}(\la_1^2,\dots,\la_n^2)$, we have $\op{deg}t^*_{\mu}=2\>|\mu|$,
so that $\op{deg}\Bbb T_{\mu}=2\>|\mu|$.

Now it is clear that Theorems 2.2--2.3 are equivalent to the following
claims about the polynomials $s^*_{\mu}$ and $t^*_{\mu}$.

\proclaim{Theorem 2.5} The polynomials $s^*_{\mu}(\la_1,\dots,\la_n)$ and
$t^*_{\mu}(\la_1,\dots,\la_n)$ form a basis of the algebras
$\La^*(n)$ and $\roman{M}^*(n)$, respectively. They can be characterized\/
\rom(within a scalar multiple\/\rom) by the following properties\/\rom:
\roster
\item"(i)" $s^*_{\mu}\in\La^*(n)$
and $\op{deg}s^*_{\mu}=|\mu|$, $t^*_{\mu}\in
\roman{M}^*(n)$ and $\op{deg}t^*_{\mu}=2\>|\mu|$\rom;
\item"(ii)" $s^*_{\mu}(\la)=0$ and
$t^*_{\mu}(\la)=0$ for all partitions $\la$ such that $l(\la)\le n$,
$|\la|\le|\mu|$, $\la\ne\mu$\rom;
\item"(iii)" $s^*_{\mu}(\mu)\ne0$ and
$t^*_{\mu}(\mu)\ne0$.
\endroster
\endproclaim

\demo{Proof} The part of these claims concerning the polynomials
$s^*_{\mu}$ was established in \cite{Ok1}; see also \cite{OO1,
Theorem 3.3}. So we shall only consider the case of $t^*_{\mu}$
(which, however, is quite similar to that of $s^*_{\mu}$).

Let $\roman{M}(n)$ denote the algebra of polynomials in $n$ variables
$\la_1,\dots,\la_n$, invariant under the natural action of the
hyperoctahedral group $H(n)$ (the generators of $\{\pm1\}^n\subset
H(n)$ change the sign of the variables). Clearly, the polynomials
$$
t_{\mu}(\la_1,\dots,\la_n):=s_{\mu}(\la_1^2,\dots,\la_n^2)
$$
form a homogeneous basis in $\roman{M}(n)$.

On the other hand, $\roman{M}(n)$ is canonically isomorphic to the graded
algebra associated with the filtered algebra $\roman{M}^*(n)$. Since
$t_{\mu}$ coincides with the top component of $t^*_{\mu}$, we
conclude that the polynomials $t^*_{\mu}$ form a basis in $\roman{M}^*(n)$.

Let us check that the polynomials $t^*_{\mu}$ satisfy the conditions
(i)--(iii). The condition (i) is evident and the conditions (ii),
(iii) will be deduced from the explicit formula
$$
\split
t^*_{\mu}(\la_1,\dots,\la_n)
&=s_{\mu}(l_1^2,\dots,l_n^2\mi\e^2,(\e+1)^2,\dots)\\
&=\frac{\op{det}\>[(l_i^2\mi\e^2,(\e+1)^2,\dots)^{\mu_j+n-j}]}{\prod_{i<j}
(l_i^2-l_j^2)}\.
\endsplit
$$
Note that the denominator is nonzero, because $l_1^2,\dots,l_n^2$
strictly decrease. Thus, we have to analyze the vanishing properties
of the numerator.

Assume $|\la|\le|\mu|$ and $\la\ne\mu$. Then there exists
$k\in\{1,\dots,n\}$ such that $\la_k<\mu_k$, whence
$$
\la_i<\mu_j\quad\text{for }\,j\le k\le i\.
$$
For such couples $i$, $j$ the $(i,j)$\<th entry of the determinant in
the numerator is
$$
(l_i^2\mi\e^2,(\e+1)^2,\dots)^{\mu_j+n-j}=\prod_{r=1}^{\mu_j+n-j}
((\la_i+n-i+\e)^2-(\e+r-1)^2)=0,
$$
so that the determinant vanishes.

If $\la=\mu$ then a similar argument shows that the matrix in the
numerator is strictly upper triangular with nonzero diagonal entries,
so that the determinant is nonzero.

Thus, we have verified the properties (i)--(iii) and it remains to
prove that polynomials with such properties are unique. For
$d=2,4,6,\dots$ consider the subspace 
$\roman{M}^*_d(n)\subset \roman{M}^*(n)$ formed
by elements of degree $\le d$. The existence property verified above
implies that the linear functionals on $\roman{M}^*_d(n)$ of the form
$$
f\mapsto f(\la),\qquad|\la|\le d,
$$
are linearly independent. On the other hand, their number just equals
$\op{dim}\roman{M}^*_d(n)$, which implies the uniqueness claim.\qed
\enddemo

Note that Theorem 2.5 also can be directly deduced from the binomial
formula, cf. second proof of Theorem 5.1 in \cite{OO1}.

\head\S3.\enspace Coherence property\endhead

Until now we fixed a classical Lie algebra $\g$. In this section
we deal with the whole series $\{\g(n)\}$ of classical Lie
algebras (of type $A,C,B$ or $D$). We define natural embeddings
$Z(\g(n))\to Z(\g(n+1))$ and we show that the elements of the
canonical bases are stable (within scalar multiples) with respect to
these embeddings.

Since the Lie algebra will vary (inside a classical series), we shall
use more detailed notation for the polynomials in $n$ variables
introduced above:
$$
s_{\mu|n},\;s^*_{\mu|n},\;t_{\mu|n},\;t^*_{\mu|n},
$$
and the elements of the distinguished basis of $Z(\g(n))$ will be
denoted by $\Bbb S_{\mu|n}$ or $\Bbb T_{\mu|n}$.

For any $n$, we shall identify $\C^n$ with the subspace of $\Bbb
C^{n+1}$ spanned by first $n$ basis vectors. Since $\g(n)$ is
realized in $\C^n$ (for series A) or in $\C^n\oplus\Bbb
C^n$ (for series C, B, D), we obtain a natural Lie algebra embedding
$\g(n)\to\g(n+1)$, which induces a natural algebra
embedding $U(\g(n))\to U(\g(n+1))$.

Next, we define the {\it averaging operator}
$$
\op{Avr}_{n,n+1}\:Z(\g(n))\to Z(\g(n+1))
$$
as the composition
$$
Z(\g(n))\hookrightarrow U(\g(n))\hookrightarrow U(\g(n+1))\to Z(\g(n+1)),
$$
where the latter arrow is the projection $\#$ defined in (2.1).
Similarly, for any couple $n<m$, we define the averaging operator
$$
\op{Avr}_{nm}\:Z(\g(n))\to Z(\g(m)),
$$
which also coincides with composition of averagings
$$
Z(\g(n))\to Z(\g(n+1))\to\dots\to Z(\g(m))\.
$$

\proclaim{Theorem 3.1} Let $c(n,\mu)$ and $c_{\pm}(n,\mu)$ be the
normalizing factors occurring in the binomial formula\/ 
\rom(see\/ \rom{(1.7)},
\rom{(1.9))}. For any $n<m$ we have
$$
\align
\op{Avr}_{nm}\bigg(\frac1{c(n,\mu)}\,\Bbb S_{\mu|n}\bigg)
&=\frac1{c(m,\mu)}\,\Bbb S_{\mu|m},\\
\op{Avr}_{nm}\bigg(\frac1{c_{\pm}(n,\mu)}\,\Bbb T_{\mu|n}\bigg)
&=\frac1{c_{\pm}(m,\mu)}\,\Bbb T_{\mu|m}\.\endalign
$$
\endproclaim

We call this the {\it coherence property\/} of the basis $\{\Bbb
S_{\mu}\}$ or $\{\Bbb T_{\mu}\}$. For the series A the coherence
property was established in our paper \cite{OO1, \S10}, by three
different methods which can 
also be transferred to the series C, B, D.
We shall give now one of the proofs, which is based on the binomial
formula.

\demo{Proof} We assume $\{\g(n)\}$ is one the series C, B, D.
First of all we shall restate Theorem 1.2. Let $g\in G$ and
$X\in\g$ be related by
$$
g^{1/2}-g^{-1/2}=X,
$$
where $g$ is close to $e\in G$ and $X$ is close to $0\in\g$.
Let $f\in I(G)$ be a test function and let $f'=f'(X)$ be the
corresponding function in a neighborhood of $0\in\g$. Write $\pm
x_1,\pm x_2,\dots$ for the eigenvalues of the matrix $X$ and define
$$
T_{\mu}(X)=s_{\mu}(x_1^2,x_2^2,\dots),\qquad X\in\g,\;
l(\mu)\le\op{rank}\g\.
$$
One can expand $f'(X)$ into a series of the invariant polynomials
$T_{\mu}(X)$, which converges in a neighborhood of the origin. We
claim that this expansion can be written as
$$
f'(x)=\sum_{l(\mu)\le\op{rank}\g}\lan\wt{\Bbb{T}}_{\mu},f\ran\>
T_{\mu}(X),\tag3.1
$$
where
$$
\wt{\Bbb{T}}_{\mu}:=\frac1{c_{\pm}(\op{rank}\g,\mu)}\,\Bbb T_{\mu}\.
$$

Indeed, since $I(G)$ is spanned by the characters $\chi_{\la}$ (here
we use the same abbreviation as in Theorem 1.2), we may assume
$f=\chi_{\la}/\chi_{\la}(e)$. But then our claim is just a
restatement of the binomial formula. To see this we have to recall
the definition of the elements $\Bbb T_{\mu}$ and use the fact that
for any element $A\in Z(\g)$ (in particular, for $A=\Bbb
T_{\mu}$), $\lan A,\chi/\chi(e)\ran$ equals the eigenvalue of $A$ in
the $U(\g)$-module corresponding to $\chi_{\la}$.

Now let $n<m$ and let $G(n)\hookrightarrow G(m)$ be the natural
embedding. Let $f\in I(G(m))$ be a test function and $f\mi G(n)$ be its
restriction to $G(n)$, which is clearly an element of $I(G(n))$. Let
us write the expansion (3.1) for $f$ and $f\mi G(n)$:
$$
\alignat2
f'(X)&=\sum_{l(\mu)\le m}\lan\wt{\Bbb{T}}_{\mu|m},f\ran\>T_{\mu|m}(X),
&\qquad X&\in\g(m),\\
(f\mi G(n))'(Y)&=\sum_{l(\mu)\le n}\lan\wt{\Bbb{T}}_{\mu\mi n},f\mi G(n)
\ran\>T_{\mu|n}(Y),&\qquad Y&\in\g(n)\.
\endalignat
$$

Assume $X=Y\in\g(n)$. Then the left-hand sides of both expansions
coincide. On the other hand, by the definition of the polynomial
functions $T_{\mu}$ and the stability property of the Schur
functions, we have
$$
T_{\mu|m}(Y)=\cases T_{\mu|n}(Y),&\text{if }\,l(\mu)\le n,\\
0,&\text{if }\,l(\mu)>n\.\endcases
$$
It follows
$$
\lan\wt{\Bbb{T}}_{\mu|n},f\mi G(n)\ran
=\lan\wt{\Bbb{T}}_{\mu|m},f\ran,\qquad f\in I(G(m)),\;l(\mu)\le n\.
$$

More generally, for any $f\in\C\>[G]$ we have
$$
\lan\wt{\Bbb{T}}_{\mu|m},f\ran=\lan\wt{\Bbb{T}}_{\mu|m},f^\#\ran
=\lan\wt{\Bbb{T}}_{\mu|n},f^\#\mi G(n)\ran,
$$
which is equivalent to the claim of the theorem.\qed
\enddemo

\remark{Remark\/ \rm3.2} It was shown in \cite{OO1},
that the coherence property
of the basis $\{\Bbb S_{\mu}\}$ is equivalent to the following
relation satisfied by the polynomials $s^*_{\mu}$: for any partition
$\mu$ of length $\le n$ and any signature $\La$ of length $n+1$,
$$
\frac{\chi_{\La}^{gl(n+1)}(e)}{c(n+1,\mu)}\,s^*_{\mu|{n+1}}(\La_1,\dots,
\La_{n+1})=\sum_{\la\prec\La}\frac{\chi_{\la}^{gl(n)}(e)}
{c(n,\mu)}\,s^*_{\mu|n}(\la_1,\dots,\la_n),\tag3.2
$$
where $\la$ denote a signature of length $n$ and $\la\prec\La$ means
the Gelfand--Tsetlin interlacing condition
$$
\La_1\ge\la_1\ge\La_2\ge\cdots\ge\La_n\ge\la_n\ge\La_{n+1}\.
$$
A similar relation, which is equivalent to the coherence property of
$\{\Bbb T_{\mu}\}$, holds for the polynomials $t^*_{\mu}$:
$$
\frac{\chi_{\La}(e)}{c_{\pm}(n+1,\mu)}\,t^*_{\mu|n+1}(\La_1,\dots,
\La_{n+1})=\sum_{\la}[\chi_{\La}:\chi_{\la}]\,\frac{\chi_{\la}(e)}
{c_{\pm}(n,\mu)}\,t^*_{\mu|n}(\la_1,\dots,\la_n)\.\tag3.3
$$
Here $\La$ is a partition of length $\le n+1$, $\mu$ is a partition
of length $\le n$, $\chi_{\La}$ is one of the characters of $G(n+1)$
(i.e., $\chi_{\La}^{sp(2n+2)}$, $\chi_{\La}^{so(2n+3)}$ or
$\chi_{\La}^{o(2n+2)}$), $\chi_{\la}$ has the same meaning (but for
the group $G(n)$ of rank $n$), and $[\chi_{\La}:\chi_{\la}]$ denotes
the multiplicity of $\chi_{\la}$ in the decomposition of $\chi_{\La}$
as restricted to $G(n)\subset G(n+1)$ \footnote{This multiplicity is given by
the well-known `branching rules' for symplectic and orthogonal
groups, see \v{Z}elobenko \cite{Z, \S\S 129--130}.}.
\endremark

In \cite{OO1, \S 10} we gave a direct derivation of the relation
(3.2). The relation (3.3) can be directly verified by a similar (but more
complicated) argument. This gives another approach to the coherence
property.

\head\S4.\enspace A distinguished basis in $I(\g)$\endhead

In this section, we study a distinguished basis in $I(\g)\subset
S(\g)$ (see section $0$ for the definition of $I(\g)$).

We equip $\g$ with an invariant inner product: for any
$X,Y\in\g$,
$$
\lan X,Y\ran=\cases\op{tr}XY,&\text{for the series A},\\
\frac12\op{tr}XY,&\text{for the series C, B, D}\.\endcases
$$
We shall identity $\g$ with its dual space $\g^*$ by making
use of the product $\lan\bcdot,\bcdot\ran$. This will allow us to interpret
$S(\g)$ as the algebra of polynomial functions on $\g$ and
then elements of $I(\g)$ will become invariant polynomial
functions on $\g$.

Assume $X$ ranges over $\g$. In the case of the series A we
shall denote by $x_1,\dots,x_n$ the eigenvalues of $X$.
In the case of the series C, D the eigenvalues of $X$ may be
written as $\pm x_1,\dots,\pm x_n$, and for the series B one has to
add one~$0$. In this notation we set
$$
\align
S_{\mu}(X)&=S_{\mu|n}(X)=s_{\mu}(x_1,\dots,x_n),\\
T_{\mu}(X)&=T_{\mu|n}(X)=t_{\mu}(x_1,\dots,x_n)=s_{\mu}(x_1^2,\dots,x_n^2)\.
\endalign
$$
Clearly, $\{S_{\mu}\}$ or $\{T_{\mu}\}$ is a homogeneous basis in
$I(\g)$.

Note that $I(\g)$ may be identified with
$\op{gr}Z(\g)$, the graded algebra associated with the fibered
algebra $Z(\g)$. The next claim is immediate from definitions.

\proclaim{Proposition 4.1} Under the identification
$I(\g)=\op{gr}Z(\g)$, each basis element 
$S_{\mu}\in I(\g)$ or $T_{\mu}\in I(\g)$ coincides with the
leading term of the basis element $\Bbb S_{\mu}\in Z(\g)$ or
$\Bbb T_{\mu}\in Z(\g)$, respectively. \qed
\endproclaim

Now write $\g(n)$ instead of $\g$ and consider the chain
$\cdots\hookrightarrow\g(n)\hookrightarrow\g(n+1)
\hookrightarrow\cdots$ of natural embeddings that we already
discussed in \S3. Note that these embedding are isometric with
respect to the inner product $\lan\bcdot,\bcdot\ran$.
Therefore, for each couple
$n<m$ there is a natural projection $\g(m)\to\g(n)$
and so an algebra morphism $S(\g(m))\to S(\g(n))$, which
in turn induces an algebra morphism
$$
\op{Proj}_{nm}\:I(\g(m))\to I(\g(n)),\qquad n<m\.
$$

\proclaim{Proposition 4.2 \rm(Stability of the bases $\{S_{\mu}\}$,
$\{T_{\mu}\}$)} Assume $n<m$. If $l(\mu)\le n$ then
$$
\op{Proj}_{nm}(S_{\mu|n})=S_{\mu|n},\qquad
\op{Proj}_{nm}(T_{\mu|m})=T_{\mu|n}\.
$$
If $n<l(\mu)\le m$ then the result is zero.
\endproclaim

\demo{Proof} This follows at once from the stability property of the
Schur polynomials.\qed
\enddemo

On the other hand, for any $n<m$ we can define linear maps
$$
\op{Avr}_{nm}\:I(\g(n))\to I(\g(m))
$$
(the {\it averaging operators\/}) just in the same way as for the
invariants in the enveloping algebras, see \S3.

\proclaim{Proposition 4.3 \rm(Coherence property)} Assume $n<m$ and
$l(\mu)\le n$. Then
$$
\op{Avr}_{nm}(S_{\mu|n})=\frac{c(n,\mu)}{c(m,\mu)}\,S_{\mu|m},\qquad
\op{Avr}_{nm}(T_{\mu|n})=\frac{c_{\pm}(n,\mu)}{c_{\pm}(m,\mu)}\,T_{\mu|m}\.
$$
In particular, the averaging operators are injective.
\endproclaim

\demo{Proof} This follows from Theorem 3.1 and Proposition 4.1.
\enddemo

\proclaim{Theorem 4.4 \rm(Characterization of the basis)} Let $n$ be fixed
and $\mu$ range over partitions of length $\le n$. The basis elements
$S_{\mu|n}$ \rom(for the series A\/\rom) or $T_{\mu|n}$ \rom(for the series
C, B, D\/\rom) are the unique, within a scalar multiple, elements of
$I(\g(n))$ that are eigenvectors of all linear mappings
$$
\op{Proj}_{nm}\circ \op{Avr}_{nm}\:I(\g(n))\to I(\g(n)),\qquad
m=n+1,n+2,\dots\.
$$
\endproclaim

\demo{Proof} By Propositions 4.2 and 4.3, we have for any $m>n$
$$
\align
\op{Proj}_{nm}(\op{Avr}_{nm}(S_{\mu|n}))
&=\frac{c(n,\mu)}{c(m,\mu)}\,S_{\mu|n},\\
\op{Proj}_{nm}(\op{Avr}_{nm}(T_{\mu|n}))
&=\frac{c_{\pm}(n,\mu)}{c_{\pm}(m,\mu)}\,T_{\mu|n}\.
\endalign
$$
Thus, each basis element is an eigenvector for all the maps
$\op{Proj}_{nm}\circ \op{Avr}_{nm}$. It remains to show that the eigenvalues
corresponding to two different basis vectors are distinct for certain
$m$. But this follows from the explicit expressions for $c(n,\mu)$
and $c_{\pm}(n,\mu)$ given in \S1. \qed
\enddemo

In \cite{OO1, (2.10)}, we obtained an explicit expression of the basis
elements $S_{\mu|n}\in I(\frak{gl}(n))$ through the matrix limits
$E_{ij}\in\frak{gl}(n)$, $1\le i,j\le n$:

\proclaim{Proposition 4.5} Let $\mu$ be a partition of length $\le n$,
$k=|\mu|$, and $\chi^{\mu}$ be the irreducible character of the
symmetric group $\goth{S}(k)$ that is indexed by $\mu$. Then
$$
S_{\mu|n}=(k!)^{-1}\sum_{i_1,\dots,i_k=1}^n\sum_{s\in\goth{S}(k)}
\chi^{\mu}(s\)E_{i_1i_{s(1)}}\cdots E_{i_ki_{s(k)}}\.
$$
\endproclaim

We shall give now a similar formula for the series C, B, D. Assume
$\g=\g(n)$ is of type C, B, D and realize it as an
involutive subalgebra in $\frak{gl}(N,\C)$, where $N=2n$ or
$N=2n+1$, as indicated in \S1. We shall assume the canonical
basis in $\C^N$ is labeled by the numbers $i=-n$,
$-n+1,\dots,n-1,n$, where $0$ is included for the series B only, and
we put
$$
F_{ij}=E_{ij}-\th_{ij}E_{-j,-i},
$$
where $\th_{ij}\equiv1$ for the series B, D and
$\th_{ij}=\op{sgn}\(i\)\op{sgn}\(j)$ for the series C. Note that for our
embedding $\g(n)\hookrightarrow\frak{gl}(N,\C)$ the elements
$F_{ij}$ form a basis in $\g(n)$.

Given a partition $\mu\vdash k$ we define a central function $\varphi^\mu$
on the symmetric group $\goth{S}(2k)$ as follows. If a permutation
$s\in\goth{S}(2k)$ is such that all its cycles are even (so that the
circle type of $s$ may be written as $2\rho$ for a certain $\rho\vdash k$)
then we put
$$
\varphi^\mu(s)=\chi^{\mu}_{\rho}
$$
where $\chi_s^{\mu}$ denotes the value of the irreducible character
$\chi^{\mu}$ at any permutation with cycle type $\rho$. Otherwise
(i.e., if $s\in\goth{S}(2k)$ has at least one cycle of odd length)
$\varphi^\mu(s)=0$.

The next claim is an exact analog of Proposition 4.5.

\proclaim{Proposition 4.6} Assume $\g=\g(n)$ is of type $C,B,D$.
Let $\mu$ be an arbitrary partition of length $\le n$ and $k=|\mu|$.
Then
$$
T_{\mu|n}=\frac1{(2k)!}\sum_{i_1,\dots,i_{2k}=-n}^n
\sum_{s\in\goth{S}(2k)}\varphi^\mu(s\)
F_{i_1i_{s(1)}}\cdots F_{i_{2k}i_{s(2k)}}\.
$$
\endproclaim

\demo{Proof} Let us view $T_{\mu|n}$ as an invariant polynomial
function $T_{\mu|n}(X)$, where $X$ ranges over $\g(n)$ and denote
by $\pm x_1,\dots,\pm x_n$ the eigenvalues of $X$ (for the series B
we exclude the zero eigenvalue). We shall use standard notation of
Macdonald's book \cite{M2}: $p_\rho$ are the power sums and
$$
z_\rho=1^{k_1}k_1!\,2^{k_2}k_2!\,\cdots\quad\text{for }\,
\rho=(1^{k_1}2^{k_2}\cdots)\.
$$

By definition of the basis $\{T_{\mu}\}$ we have
$$
\align
T_{\mu|n}(X)=s_{\mu}(x_1^2,\dots,x_n^2)&=\sum_{\rho\vdash k}z_{\rho}^{-1}
\chi_{\rho}^{\mu}p_{\rho}(x_1^2,\dots,x_n^2)\\
&=\sum_{\rho\vdash k}z_{\rho}^{-1}2^{-l(\rho)}\chi_{\rho}^{\mu}p_{2\rho}
(x_1,-x_1,\dots,x_n,-x_n)\\
&=\sum_{\rho\vdash k}z_{2\rho}^{-1}\chi_{\rho}^{\mu}p_{2\rho}(x_1,-x_1,\dots,
x_n,-x_n)\.
\endalign
$$

Let $s$ range over elements of $\goth{S}(2k)$ with even cycle type $2\rho$
(where $\rho\vdash k$ depends on $s$). Then the latter expression can be
written as
$$
T_{\mu|n}(X)=\frac1{(2k)!}\sum_s\varphi^\mu(s\)
p_{2\rho}(x_1,-x_1,\dots,x_n,-x_n)\.
$$

On the other hand, for any (even) $m$, the invariant polynomial
function
$$
p_m(X)=\op{tr}(X^m)=p_m(x_1,-x_1,\dots,x_n,-x_n)
$$
corresponds to the element
$$
\sum_{i_1,\dots,i_m=-n}^nF_{i_1i_2}F_{i_2i_3}\cdots F_{i_mi_1}\in I(\g(n)),
$$
which also can be written as
$$
\sum_{i_1,\dots,i_m=-n}^nF_{i_1i_{s(1)}}\cdots F_{i_mi_{s(m)}}
$$
for any cyclic permutation $s$ of the indices $1,\dots,m$. This
concludes the proof.\qed
\enddemo

There exists another characterization of the bases $\{S_{\mu}\}$,
$\{T_{\mu}\}$, which is based on the following observation, which is
undoubtedly well-known.

\proclaim{Proposition 4.7} Define a classical group $\wt G$ containing $G$
as follows.

If $G=GL(n,\C)$ then $\wt G=GL(n,\C)\times GL(n,\C)$ and
the embedding $G\hookrightarrow\wt G$ has the form
$g\mapsto(g,(g')^{-1})$.

If $G$ is one of the groups $SO(2n+1,\C)$ or $Sp(2n,\C)$,
$SO(2n,\C)$ then $\wt G=GL(N,\C)$, where $N=2n+1$ or $N=2n$,
respectively.

Then the adjoint action of $G$ in $\g$ can be extended to a
linear action of the group $\wt G\supset G$ such that the induced
representation of $\wt G$ in the symmetric algebra $S(\g)$ is a
multiplicity free polynomial representation.
\endproclaim

\demo{Proof} Assume $G=GL(n,\C)$. Using the bijection
$E_{ij}\leftrightarrow e_j\otimes e_i$, where $\{e_i\}$ stands for
the natural basis of $\C^n$, we identify $\g=\goth{gl}(n,\C)$
with $\C^n\otimes\C^n$. The action of the group $\wt
G=GL(n,\C)\times GL(n,\C)$ in the vector space $\g=\C^n\otimes\C^n$
is the natural one. Clearly, its restriction to
the subgroup $G=\{(g,(g')^{-1})\}$ is equivalent to the adjoint
action. It is well-known that action of the group $GL(n,\C)\times
GL(n,\C)$ in $S(\C^n\otimes\C^n)$ is a multiplicity free
polynomial representation, see, e.g., \v{Z}elobenko \cite{Z,\S56} or Howe
\cite{H, 2.1}.

Assume $G=Sp(2n,\C)$ and put $N=2n$. One can identify
$\g=sp(2n,\C)$ with $S^2(\C^n)$ in such a way that the adjoint
action of $G$ in $\g$ will coincide with the natural action of
$\wt G=GL(N,\C)$ in $S^2(\C^N)$, restricted to $G\subset\wt
G$. It is well known \cite{H, 3.1} 
that the natural action of $GL(N,\C)$ in
$S(S^2(\C^N))$ is a multiplicity free polynomial representation.

Assume now $G=SO(N,\C)$, where $N=2n+1$ or $N=2n$. Then one can
identity $\g$ with $\La^2(\C^N)$ in such a way that the
adjoint action of $G$ in $\g$ will coincide with the restriction
to $G$ of the natural action of the group $\wt G=GL(N,\C)$ in
$\La^2(\C^N)$. On the other hand, it is well known \cite{H, 3.8}
that $S(\La^2\C^N)$ is a multiplicity free polynomial $GL(N,\Bbb
C)$-module.

Finally, note that in this argument, we could replace the special
orthogonal group $SO(N,\C)$ by the complete orthogonal group
$O(N,\C)$.\qed
\enddemo

The next result provides us with an alternative characterization of
the bases $\{S_{\mu}\}$, $\{T_{\mu}\}$.

\proclaim{Theorem 4.8} Let $G$ be a classical group and $\wt G\supset G$ be
as in Proposition\/ \rom{4.7}. The basis elements $S_{\mu}\in I(\g)$ or
$T_{\mu}\in I(\g)$ are the unique\/ \rom(within a scalar multiple\/\rom)
elements of $I(\g)\subset S(\g)$ that generate, under the
action of $\wt G$, irreducible submodules of $S(\g)$.
\endproclaim

\demo{Proof} Examine the case $G=GL(n,\C)$. The representation of
the group $\wt G=GL(n,\C)\times GL(n,\C)$ in the space
$S(\C^n\otimes\C^n)$ is decomposed in the direct sum of the
irreducibles $V_{\mu|n}\otimes V_{\mu|n}$, where $\mu$ ranges over
partitions of length $\le n$ (see \cite{Z, \S56} or \cite{H, 2.1.2}).
Each component $V_{\mu|n}\otimes
V_{\mu|n}$ clearly contains a unique (within a scalar multiple)
vector $S'_{\mu}$, invariant under the subgroup $G=\{(g,(g')^{-1})\}$.
It is also evident that the elements $S'_{\mu}$ can be characterized
as unique elements of $I(\g)$ generating irreducible
$G$-submodules of $S(\g)$. Thus, we have to check that
$S'_{\mu}=S_{\mu}$ (within a scalar multiple).

By Theorem 4.4, it suffices to show that each $S'_{\mu}$ is an
eigenvector for the maps $\op{Proj}_{nm}\circ \op{Avr}_{nm}$,
$m=n+1,n+2,\dots$,
but this follows from the fact that
$$
\op{Avr}_{nm}(V_{\mu|n})\subset V_{\mu|m},
\qquad \op{Proj}_{nm}(V_{\mu|m})=V_{\mu|n},
$$
which in turn follows from the very definition of these maps.

For other groups $G$ the argument is similar.

Let $G=Sp(2n,\C)$ and put $N=2n$. The representation of $\wt
G=GL(N,\C)$ in $S(S^2(\C^N))$ is the direct sum of the
irreducibles $V_{M|N}$, where $M$ is a Young diagram with even rows (the
number of rows does not exceed $N$); on the other hand, an
irreducible polynomial $GL(N,\C)$-module $V_{M|N}$ contains a
vector invariant under $Sp(2n,\C)$, where $n=N/2$, if and only if
all the columns of $M$ are even, and then such a vector is unique
within a scalar multiple (see \cite{H, 3.1} or \cite{M2, VII, (6.11)}).

Thus, an irreducible component $V_{M|N}\subset S(S^2(\C^{2n}))$
has a nonzero (and then one-dimensional) intersection with $I(\g)$
if and only if $M$ can be written as
$$
M=2\mu\cup2\mu=(2\mu_1,2\mu_1,2\mu_2,2\mu_2,\dots,2\mu_n,2\mu_n),
$$
where $\mu$ is a partition of length $\le n$. Let $T'_{\mu}$ be any
nonzero element of the one-dimensional space $V_{2\mu\cup2\mu|2n}\cap
I(\g)$. The same argument as above shows that $T'_{\mu}=T_{\mu}$,
within a scalar multiple.

(Note, however, a difference with the case $G=GL(n,\C)$ examined
above. For $G=GL(n,\C)$, we saw that each irreducible $\wt
G$-submodule contained a $G$-invariant, while now only a part of
components possess $G$-invariants.)

Now let $G=SO(N,\C)$, where $N=2n+1$ or $N=2n$, and remark that
in both cases $I(\g)$ coincides with the subspace of $S(\g)$
formed by the $G'$-invariants, where $G'=O(N,\C)$. The
representation of $\wt G=GL(N,\C)$ in the space $S(\La^2(\Bbb
C^N))$ decomposes into the direct sum of irreducible modules
$V_{M|N}$, where $M$ is a diagram with even columns. On the other hand,
a $G'$-invariant in $V_{M|N}$ exists (and then is unique, within a
scalar multiple) if and only if $M$ has even rows 
(see \cite{H, 3.8} or \cite{M2, VII, (3.14)}). Thus we again
obtain that $M$ has the form $2\mu\cup2\mu$ and conclude the proof as
above.\qed
\enddemo

\head\S5.\enspace A relationship between the bases in $Z(\g)$
and $I(\g)$\endhead

Here we exhibit a linear isomorphism $S(\g)\to U(\g)$,
called the special symmetrization, which maps $I(\g)$ onto
$Z(\g)$ and takes the canonical basis in $I(\g)$ to that
of $Z(\g)$. The main result is Theorem 5.2. We reduce it to
Proposition 5.3 which in turn is reduced to Proposition 5.4. These
two propositions are of independent interest

We shall heed the concept of generalized symmetrization proposed by
Olshanski \cite{O2}. Let us identify the algebra $U(\g)$ with
a subspace of the dual of $\Cal O_e(G)$ (see the beginning of \S2).
Similarly, we shall identify $S(\g)$ with a 
subspace of the dual space to $\Cal O_0(\g)$, the space of germs of
holomorphic functions defined at a neighborhood of the origin
$0\in\g$.

Assume we are given a map $F\:\g\to G$ with the following
properties:
\roster
\item"(i)" $F$ is holomorphic and defined in a neighborhood of the origin,
invariant under the adjoint representation;
\item"(ii)" $F$ takes $0\in\g$ to $e\in G$;
\item"(iii)" the differential of $F$ at the origin is the
identical map $\g\to\g$;
\item"(iv)" $F$ is equivariant with
respect to the group $G$ acting by the adjoint representation in
$\g$ and by conjugations on itself.
\endroster

Then $F$ induces an isomorphism of local rings $\Cal O_0(\g)\to
\Cal O_e(G)$, which by duality determines a linear isomorphism
$$
\si\:S(\g)\to U(\g)\.
$$

Note that this construction has a sense for any complex Lie group; in
fact it also works for formal groups.

When $F$ is the exponential map, the corresponding map $\si$
coincides with the standard symmetrization map. Following \cite{O2},
in the general case we call $\si$ a {\it generalized symmetrization}.

A generalized symmetrization $\si$ shares a number of properties of
the standard symmetrization. In particular, $\si$ preserves leading
terms and also induces a bijection between invariants of the group.

Thus, for a classical group $G$, $\si$ induces a linear isomorphism
$I(\g)\to Z(\g)$ (for the series D we shall assume $F$ is
invariant under $G'\supset G$).

\proclaim{Theorem 5.1 \rm(\cite{OO1, Theorem 14.1})} Assume 
$\g=\frak{gl}(n,\C)$ and $F(X)=1+X$, and let $\si$ be the corresponding
generalized symmetrization. Then
$$
\si(S_{\mu|n})=\Bbb S_{\mu|n}\quad
\text{for any }\,\mu,\;\ell(\mu)\le n\.\qed
$$
\endproclaim

This generalized symmetrization was introduced in \cite{O1}; then it
was used in the note \cite{KO} and called the {\it special
symmetrization}. We aim to find a similar generalized symmetrization
for the series C, B, D.

\proclaim{Theorem 5.2} Let $\g=\g(n)$ be any classical Lie
algebra of type $C,B,D$, realized as an involutive subalgebra of
$\frak{gl}(N,\C)$, where $N=2n$ or $N=2n+1$. Define $F\:\g\to G$
by the relation
$$
F(X)^{1/2}-F(X)^{-1/2}=X\quad\text{for }\,X\in\g,
$$
i.e.,
$$
\align
F(X)&=1+X^2\!/2+((1+X^2\!/2)^2-1)^{1/2}
=1+X^2\!/2+X(1+X^2\!/4)^{1/2}\\
&=1+X+\cdots,
\endalign
$$
where the eigenvalues of the matrix $X$ are supposed to be
sufficiently small.

Let $\si$ be the corresponding generalized symmetrization. Then
$$
\si(T_{\mu|n})=\Bbb T_{\mu|n}\quad\text{for any }\,\mu,\;l(\mu)\le n\.
$$
\endproclaim

Note that $F$ is indeed a well-defined (local) map from $\g$ to
$G$. We call $\si$ the {\it special symmetrization\/} for the
classical Lie algebras of type C, B, D. We refer to \cite{O2} for
explicit formulas for $\si\:S(\g)\to U(\g)$ and its inverse
$\si^{-1}\:U(\g)\to S(\g)$.

Note that the proof of the theorem given below differs from the
argument used in \cite{OO1} for the series A. On the other hand, the
argument given below also holds for the series A. We do not know if
the approach of \cite{OO1, Theorem 14.1}, can be transferred to the
series C, B, D.

\demo{Proof} Let $\varphi\in \Cal O_e(G)$ stand for a test element. By
definition, we have to prove the equality
$$
\lan\Bbb T_{\mu|n},\varphi\ran=\lan T_{\mu|n},\varphi\circ F\ran,
$$
where the brackets on the left-hand side denote the pairing between $U(\g)$
and $\Cal O_e(G)$ while the brackets on the right-hand side denote the pairing
between $S(\g)$ and $\Cal O_0(\g)$.

Let
$$
\Cal O_e(G)^{\op{inv}}\subset \Cal O_e(G)\quad\text{and}\quad
\Cal O_0(\g)^{\op{inv}}\subset \Cal O_0(\g)
$$
denote the subspaces of invariants with respect to the action of the
group $G$ (or the group $G'$, for the series D). Choose a compact
form $K\subset G$ (or $K'\subset G'$). Averaging over $K$ (or $K'$)
determines an invariant projection
$$
\Cal O_e(G)\to \Cal O_e(G)^{\op{inv}}.
$$
Using it we reduce our problem to the case $\varphi\in \Cal O_e(G)^{\op{inv}}$.

We shall apply the binomial formula (Theorem 1.2) which can be
conveniently written as
$$
(\varphi\circ F)(X)=\sum_{l(\mu)\le n}\frac{\lan\Bbb T_{\mu|n},\varphi\ran}
{c_{\pm}(n,\mu)}\,T_{\mu|n}(X),
$$
where $\varphi\in \Cal O_e(G)^{\op{inv}}$. Indeed, if $\varphi$ has the form
$\varphi=\chi_\la/\chi_\la(e)$, where $\chi_\la$ is one of the characters
$\chi_{\la}^{sp(2n)}$, $\chi_{\la}^{so(2n+1)}$, $\chi_{\la}^{o(2n)}$
then the above relation just coincides with the binomial formula,
because of the definition of the map $F$ and the relation
$$
\lan\Bbb T_{\mu|n},\chi_\la/\chi_\la(e)\ran=t^*_{\mu|n}(\la_1,\dots,\la_n)\.
$$
Further, since the linear span of the characters $\chi_\la$ is dense in
$\Cal O_e(G)^{\op{inv}}$ with respect to the adic topology defined by the
unique maximal ideal of $\Cal O_e(G)^{\op{inv}}$, our expansion holds for any
$\varphi\in \Cal O_e(G)^{\op{inv}}$.

On the other hand, any element $\psi\in \Cal O_0(\g)^{\op{inv}}$ can be
uniquely written as a series
$$
\psi=\sum_{l(\mu)\le n}c_{\mu}(\psi\)T_{\mu|n},\qquad c_{\mu}(\psi)\in\C,
$$
converging in the adic topology of $\Cal O_0(\g)^{\op{inv}}$. By taking
$\psi=\varphi\circ F$ and comparing the both expansion we reduce the
problem to the following relation
$$
c_{\mu}(\psi)=\frac{\lan T_{\mu|n},\psi\ran}{c_{\pm}(n,\mu)},\qquad\psi\in
\Cal O_0(\g)^{\op{inv}}.
$$

Finally, without loss of generality we may assume $\psi=T_{\nu|n}$
for a certain partition $\nu$ of length $\le n$, because
$\{T_{\nu|n}\}$ is a topological basis in $\Cal O_0(\g)^{\op{inv}}$.

Thus, the final reduction of the problem is the following claim.
\enddemo

\proclaim{Proposition 5.3} Let $\mu$, $\nu$ be partitions of length $\le n$.
Then
$$
\lan T_{\mu|n},T_{\nu|n}\ran=\dt_{\mu\nu}\,c_{\pm}(n,\mu)\.
$$
\endproclaim

\remark{Comment} The brackets in the left-hand side can be understood
in two different but equivalent ways. First, they represent the
pairing between a distribution supported at the origin and a polynomial
function, i.e.,
$$
\lan T_{\mu|n},T_{\nu|n}\ran=(\partial(T_{\mu|n}\)T_{\nu|n})(0),
$$
where $\partial(T_{\mu|n})$ stands for the differential operator on $\g$
with constant coefficients that corresponds to $T_{\mu|n}\in
S(\g)$. Second, the brackets may denote the canonical extension
to $S(\g)$ of the inner product on $\g$.
\endremark

\demo{Proof} Probably, the proposition could be proved starting from
the explicit expression for $T_{\mu|n}$ given in Proposition 4.6 but
we prefer to use another argument.

Let us extend to $S(\g)$ the inner product $\lan\bcdot,\bcdot\ran$ in $\g$.
Then one can associate with $(S(\g),\lan\bcdot,\bcdot\ran)$ a
reproducing kernel $\Cal E(X,Y)$, where $X,Y\in\g$. By
definition,
$$
\Cal E(X,Y)=\sum_{\a}\psi_{\a}(X\)\psi_{\a}^*(Y),
$$
where $\{\psi_{\a}\}$ is an arbitrary {\it homogeneous\/} basis in
$S(\g)$ and $\{\psi_{\a}^*\}$ is the dual basis. Note that $\Cal
E(X,Y)$ does not depend on the choice of the basis.

Similarly, to the inner product space $(I(\g),\lan\bcdot,\bcdot\ran)$ also
corresponds a reproducing kernel $\Cal{F}(X,Y)$.

Between these two kernels there is an evident relation
$$
\Cal{F}(X,Y)=\int\Cal E(X,\op{Ad}u\cdot Y)\,du
=\int\Cal E(\op{Ad}u\cdot X,Y)\,du,
$$
where the integral is taken over the compact form $K\subset G$ (or
$K'\subset G'$) equipped with the normalized Haar measure.

Further, taking as $\{\psi_{\a}\}$ the basis of monomials formed from
a basis in $\g$, we obtain that
$$
\Cal E(X,Y)=e^{\lan X,Y\ran}
$$
whence
$$
\Cal{F}(X,Y)=\int e^{\lan X,\op{Ad}u\cdot Y\ran}\,du\.
$$

On the other hand, the claim of the proposition means that
$$
\Cal{F}(X,Y)=\sum_{l(\mu)\le n}
\frac{T_{\mu|n}(X\)T_{\mu|n}(Y)}{c_{\pm}(n,\mu)},
$$
because $\{T_{\mu|n}\}$ is a homogeneous basis in $I(\g)$.

Thus, we have reduced our problem to the following claim.
\enddemo

\proclaim{Proposition 5.4} Let $u$
range over a compact form $K\subset G$ \rom(or
$K'\subset G'$, for the series D\/\rom) and let $du$ denote the normalized
Haar measure. The following expansion holds
$$
\int e^{\lan X,\op{Ad}u\cdot Y\ran}\,du=\sum_{l(\mu)\le n}
\frac{T_{\mu|n}(X\)T_{\mu|n}(Y)}{c_\pm(n,\mu)},\qquad X,Y\in\g=\g(n)\.\tag5.1
$$
\endproclaim

\demo{Proof} We shall show that this can be obtained from the binomial
formula (Theorem 1.2) by a limit transition.

First of all, without loss of generality we may assume $X$ and $Y$
are diagonal matrices with diagonal entries $\pm x_i$ and $\pm y_i$,
respectively (as usual, we omit the zero entry in the case of the
series D). Then in the right-hand side one can replace
$T_{\mu|n}(X)$ by $t_{\mu|n}(x_1,\dots,x_n)$ and $T_{\mu|n}(Y)$ by
$t_{\mu|n}(y_1,\dots,y_n)$. After this the right-hand side becomes
very similar to the right-hand side of the binomial formula: the only
difference is that in the binomial formula, we have
$t_{\mu|n}^*(\la_1,\dots,\la_n)$ instead of
$t_{\mu|n}(y_1,\dots,y_n)$.

Let us compare now the left-hand sides of both formulas (in the
discussion below one has to replace $K$ by $K'$ for the series D).

In the binomial formula (1.8) occurs the normalized irreducible
character $\chi_\la/\chi_\la(e)$, which can be viewed as a
spherical function for the Gelfand pair $(K\times K,K)$. The
left-hand side of formula (5.1) also can be viewed as a spherical
function for a Gelfand pair. Namely, this pair consists of the Cartan
motion group $K\ltimes\frak{k}$ (the semidirect product of $K$ and its Lie
algebra $\frak{k}$ viewed as a $K$-module) and its subgroup $K$. The vector
$x=(\pm x_i)$ can be considered as the argument of the spherical
function and $y=(\pm y_i)$ is the parameter.

Now we shall use the well-known relation between the spherical
functions of a symmetric space (say, of compact type) and the
spherical functions of the corresponding Cartan motion group (see,
e.g., Dooley and Rice \cite{DR}).

In our context, this relation looks as follows. Let $\e$ be a small
parameter which then will tend to $0$. Assume $g\in K$ has the form
$g=1+\e X+O(\e)$, where $X\in\frak{k}$, and the partition $\la$ has the
form $\la=\e^{-1}y+O(\e^{-1})$.\footnote{Since the both sides of (5.1)
are invariant with respect to the action of the hyperoctahedral group
on $y_1,\dots,y_n$, we may assume $y_1\ge\cdots\ge y_n\ge0$, so that
the approximation of the vectors $\e^{-1}y$ by partitions does
exist.} 
Then in the limit $\e\to0$ the normalized
character indexed by $\la$ turns into the spherical function for the
Cartan motion group, indexed by $y$.

Since the top homogeneous component of the polynomial $t^*_{\mu|n}$
coincides with $t_{\mu|n}$, the right-hand side of the binomial
formula (1.8) will turn, after this limit transition, into the
right-hand side of formula (5.1).

This concludes the proof of Proposition 5.4 and at the same time of
Proposition 5.3 and of Theorem 5.2.
\enddemo

\head\S6.\enspace Appendix: Bispherical functions on $Sp(2n,\C)\setminus
GL(2n,\C)/O(2n,\C)$\endhead

In this section we discuss a curious fact suggested by the results of
\S4.

Consider two Gelfand pairs $(\wt G,G)$:
$$
\align
(GL(N,\C),O(N,\C)),&\qquad N=2n+1\,\text{ or }\,N=2n,
\intertext{and}
(GL(N,\C),Sp(N,\C)),&\qquad N=2n\.
\endalign
$$
The spherical functions of these pairs are matrix elements
$g\mapsto(V(g\)\xi,\xi)$, where $g$ ranges over $\wt G$, $V$ as an
irreducible finite-dimensional $\wt G$-module, and $\xi$ ia a
$G$-invariant vectors. As is well-known these spherical functions can
be identified (with an appropriate parametrization of the double
$G$-cosets in $\wt G$) as Jack polynomials $P_{\mu}^{(\a)}$ in $n$
variables, where the parameter $\a$ takes the value $\a=2$ for the
former pair and the value $\a=1/2$ for the latter one, and $\mu$
ranges over partitions of length $\le n$.

The aim of this Appendix is to examine what happens when these two
pairs are ``mixed'' in the following sense. We set
$$
\wt G=GL(2n,\C),\quad G_1=O(2n,\C),\quad G_2=Sp(2n,\C),
$$
and assume $V$ is an irreducible $\wt G$-module admitting both a
$G_1$-invariant $\xi$ and a $G_2$-invariant $\eta$. Then the dual
module $V^*$ also possesses a $G_2$-invariant $\eta^*$ and we form
the matrix element
$$
\varphi(g)=\lan V(g\)\xi,\eta^*\ran,\qquad g\in\wt G,
$$
which we call a {\it bispherical function}. We shall calculate $\ph$
in terms of a natural parametrization of the $(G_2,G_1)$-cosets in
$\wt G$ and we shall see that $\ph$ is a Schur polynomial (note that
the Schur polynomials are the Jack polynomials with $\a=1$).

Without loss generality we may assume $V$ is a polynomial module. We
shall specify the embeddings $G_1\hookrightarrow\wt G$ as indicated
in \S1.

\proclaim{Theorem 6.1} Let $V$ be an irreducible polynomial representation
of $\wt G=GL(2n,\C)$ admitting both a $G_1$-invariant $\xi$ and a
$G_2$-invariant $\eta$ \rom(where $G_1=O(2n,\C)$
and $G_2=Sp(2n,\Bbb C)$\<\rom), i.e.,
$$
V=V_{2\mu\cup2\mu|2n},
$$
where $\mu$ is a partition of length $\le n$ and
$$
2\mu\cup2\mu=(2\mu_1,2\mu_1,\dots,2\mu_n,2\mu_n)\.
$$

Then the spherical function
$$
\ph_{\mu}(g)=\lan V_{2\mu\cup2\mu|2n}(g\)\xi,\eta^*\ran,
$$
where $g\in\wt G$ and $\eta^*$ stands for a $G_2$-invariant in $V^*$,
is proportional to the Schur polynomial $s_{\mu}$ in $n$ variables
under a suitable parametrization of the double cosets $G_2gG_1$ in
$\wt G$.
\endproclaim

\demo{First proof} Thus proof is based on Theorem 4.8. Let $g\mapsto g'$
denote transposition of $2n\times2n$ matrices that corresponds to the
symmetric form $M$ preserved by the subgroup $G_1=O(2n,\C)$, see
\S1 for the choice of $M$. I.e., $g\mapsto g'$ is transposition with
respect to the secondary diagonal. The map $g\mapsto gg'$ defines a
bijection between the cosets $gG_1\subset\wt G$ and (nondegenerate)
$2n\times2n$ matrices, symmetric relative the secondary diagonal so
one can write
$$
\ph_{\mu}(g)=\psi_{\mu}(gg'),
$$
where $\psi_{\mu}$ is a polynomial function on symmetric matrices

The group $\wt G$ acts on symmetric matrices, and under this action
$\psi_{\mu}$ is specified by two properties: first, it transforms
according to $V_{2\mu\cup2\mu}$ and, second, it is $G_2$-invariant.

Next, we identity the space of symmetric matrices with the Lie
algebra $sp(2n,\C)$ or, which is better, with its dual space.
Under this identification, the functions $\psi_{\mu}$ become elements
of the space $I(sp(2n,\C))$, and application of Theorem 4.8
implies that these are just the basis elements $T_{\mu}$ (within a
scalar multiple, as usual).

Given a $n$-tuple $(x_1,\dots,x_n)$ of nonzero complex numbers, form
a diagonal matrix,
$$
g(x)=\op{diag}\(x_1,\dots,x_n,x_n,\dots,x_1)\.
$$
Then
$$
g(x)(g(x))'=\op{diag}\(x_1^2,\dots,x_n^2,x_n^2,\dots,x_1^2)\.
$$
Under the above identification between symmetric matrices and
elements of the (dual space to) the Lie algebra $sp(2n,\C)$ the
latter matrix turns into
$$
\op{diag}\(x_1^2,\dots,x_n^2,-x_n^2,\dots,-x_1^2)
$$
which is an element of the Cartan subalgebra of $sp(2n,\C)$.

It follows
$$
\ph_{\mu}(g(x))=\psi_{\mu}(g(x)(g(x))')=\op{const}s_{\mu}(x_1^4,\dots,
x_n^4),
$$
by the definition of the elements $T_{\mu}$ for the series C. To
conclude the proof we remark that each double $(G_2,G_1)$-coset in
$\wt G$ has the form $G_2g(x)G_1$ with a certain $x=(x_1,\dots,x_n)$.\qed
\enddemo

\demo{Second proof\/ \rom(sketch\/\rom)}
One can choose a compact form $\wt{G}^u\subset\wt G$
(isomorphic to $U(2n)$) such that $G_i^u=\wt G^u\cap
G_i$ is a compact form of $G_i$ for $i=1,2$. We also can arrange so
that the matrices $g(x)$ be in $\wt G^u$ provided
$|x_1|=\cdots=|x_n|=1$. Then each $(G_2^u,G_1^u)$-coset in $\wt G^u$
is of the form $G_2^ug(x)G_1^u$, so that we can use $x_1,\dots,x_n$
as parameters of the $(G_2^u,G_1^u)$-cosets in $\wt G^u$.

A direct computation shows that the ``radial part'' of the Haar
measure of the compact group $\wt G^u$, expressed in these
parameters, has the density
$$
w(x)=\op{const}\prod_{1\le i<j\le n}|x_i^4-x_j^4|^2.
$$
It follows that the bispherical functions are symmetric orthogonal
polynomials in $x_1^4,\dots,x_n^4$ with weight $w(x)$.

Finally, analysis of the weights of $V_{2\mu\cup2\mu|2n}$ shows that
the leading term of the bispherical function (with respect to the
lexicographic order on monomials in $x_1^4,\dots,x_n^4$) is equal to
$$
\op{const}x_1^{4\mu_1}\cdots x_n^{4\mu_n}.
$$

Hence our orthogonal polynomials coincide with the Schur polynomials
$s_\mu(x_1^4,\dots,x_n^4$).\qed
\enddemo

\enddocument
\bye